\documentstyle[12pt,a4]{article}

\newcommand{\be}{\begin{equation}}
\newcommand{\ee}{\end{equation}}
\newcommand{\bea}{\begin{eqnarray}}
\newcommand{\eea}{\end{eqnarray}}
\newcommand{\unop}{1\!{\rm I}}
\date{\ }
\title{spl(p,q) superalgebra and differential operators}
\author{Y. Brihaye \\
Department of Mathematical Physics\\
University of Mons\\
Av. Maistriau, B-7000 MONS, Belgium.\\
Stefan Giller{$^{*}$},
       Piotr Kosinski \thanks{$^{\dagger}$ Work supported by grant $n^o$
KBN 2P03B07610}\\
Department of Theoretical Physics\\
University of Lodz\\
Pomorska 149/153, 90-236 Lodz, Poland}
\begin{document}
\begin{titlepage}
\maketitle
\thispagestyle{empty}
\begin{abstract}
Series of finite dimensional representations of the superalgebras spl(p,q)
can be formulated in terms of linear differential operators acting
on a suitable space of polynomials.
We sketch the general ingredients necessary 
to construct these representations and 
present examples related to spl(2,1) and spl(2,2). 
By revisiting the products of projectivised representations of sl(2), 
we are able to construct new sets of differential operators 
preserving some space of polynomials in two or more variables. 
In particular, this allows to express the representation of spl(2,1)
in terms of matrix differential operators in two variables. 
The  corresponding operators provide the building blocks 
for the construction of quasi exactly solvable systems
of two and four equations in two variables. 
We also present a quommutator deformation of spl(2,1)
which, by construction, provides an appropriate basis for analyzing 
the quasi exactly solvable systems of finite difference equations.     
\end{abstract}
\end{titlepage}
\hfill {\it In memory of our friend S. Malinowski}
\vskip 2. cm 
\section{Introduction}
The number of quantum mechanical problems for which the spectral equation can
be solved algebraically is rather limited. It is therefore not surprising
that the quasi exactly solvable (QES) equations 
\cite{tur1}, \cite{ush} attract some attention.
For these equations, indeed, a finite number of eigenvectors can be obtained
by solving an algebraic equation. The study of QES equations has motivated
the classification of finite dimensional real Lie algebras of first order
differential operators. The algebras which, in a suitable representation,
preserve a finite dimensional module of smooth real functions
are particularly relevant for QES equations. The case of
one variable was addressed and solved some years ago \cite{tur2}. 
There is, up to
an equivalence only one algebra, for instance sl(2), acting on the space of
polynomials of degree at most $n$ in the variable.
For two variables, the classification is more involved. It is described
respectively in \cite{gko1} and \cite{gko2} for complex and real variables.
The corresponding real QES operators finally 
emerge in seven classes summarized
in table 7 of \cite{gko2}.
\par
The natural next step is to classify the QES systems of two 
equations \cite{turshif,bk,bn}. 
Obviously, the number of possible finite dimensional modules of
functions to be preserved grows with the number of independent
variables admitted.
It is known that, in general, the 
underlying algebraic structure is not a (super) Lie algebra.
As far as 2$\times$2 matrix differential operators
of one variable are considered, 
the few cases where the algebra is a Lie algebra, 
it is generically sl(2)$\times$sl(2) or spl(2,1) \cite{bk}.
A generalisation of these operators to the case of V 
variables leads naturally to the symmetry algebras
sl(2)$\times$sl(V+1) or spl(V+1,1) \cite{bn}.
The first case is rather trivial while the case
where a superalgebra is involved  seems
to be worth to be studied in more details.
In the present paper we study the modules of 
polynomials is several variables 
which can serve as basis of representation to the 
superalgebras of the  spl(p,q) family \cite{hurni}.
\par
In the second section, we present some new aspects of 
the representations of the Lie algebra sl(p)
in terms of differential operators. As a byproduct we obtain a simpler
description of some of the QES operators 
classified in \cite{gko2}. 
\par
Thereafter, in section 3 we consider the superalgebra
spl(p,q) and sketch the form of its representations which are
expressible in terms of linear differential operators.
We construct explicitely some of these representations in the cases
spl(2,1) and spl(2,2).  We obtain in particular a realization
of spl(2,1) in terms of partial differential operators of two real 
variables.
This provides the algebraic basis for 
the description of certain QES systems
of two or four equations in two variables. The properties
of the relevant tensorial operators are given in section 4.
\par
The tensorial labelling used to write the generators of the 
algebras considered simplifies considerably the presentation 
of their  structure constants.
Taking advantage of this labelling,
we are able to find a deformation of spl(2,1) similar, in the
spirit, to the Witten-Woronowicz deformation of sl(2) \cite{witt,woro,zachos}.
For this deformed algebra, we also find the representations  
that we express in terms of finite difference operators.
The normal ordering rules obeyed  by these operators
are more appropriate for the classification of finite difference
QES operators than the ones given in a previous paper \cite{bgk}.
\section{Lie algebra sl(M+1)}
Consider $M$ independent real variables $\vec x \equiv x_1,...,x_M$ 
and denote $P(m,M)$
the vector space of polynomials of overal degree at most $m$ in these
variables i.e.
\be
P(m,M) = {\rm span} \lbrace x^{m_1}_1 x^{m_2}_2\cdots x^{m_M}_M\quad ;
\quad m_1+m_2 + \cdots m_M \leq m\rbrace
\ee
(for shortness we use $P(m,1)\equiv P(m)$).
The linear differential operators preserving $P(m,M)$ can be constructed
(in a sense of enveloping algebra) out of $(M+1)^2$ operators
\cite{turshif,bn}
\begin{eqnarray}
J^0_0(\vec x,m) &=& D-m+\gamma \quad , \quad 
D\equiv \sum^M_{j=1} x_j {\partial \over {\partial x_j}} \nonumber \\
J^k_0(\vec x,m) &=& {\partial \over
{\partial x_k}}\quad , \quad k=1,\ldots, M   \nonumber \\
J^0_k(\vec x,m) &=& -x_k(D-m)\quad , \quad k=1,\ldots, M \nonumber \\
J^l_k(\vec x,m) &=& - x_k{\partial\over{\partial x_l}}\quad + 
  \gamma \delta^l_k
\  , \  \quad k,l=1,\ldots, M
\label{j}
\end{eqnarray}
 which obey the commutation relations of $gl(M+1)$ :
\be
[J^{b}_{a}, J^{d}_{c}] = 
 \delta^{d}_{a} J^{b}_{c}
- \delta^{b}_{c}J^{d}_{a} 
\quad \ \ , \ \ a,b,c,d = 0,1, \ldots , M \ .
\label{comj}
\ee
irrespectively of the constant $\gamma$.
\par
For later convenience we introduce for $M=1$ the 
fundamental QES operators
 \bea
\label{sl2}
  j_+(x,m) &=& -J_1^0(x,m) = x(\partial_x-m) \nonumber \\
  j_0(x,m) &=& {1\over 2}(J_0^0(x,m)-J_1^1(x,m)) =\partial_x - {m\over2} \\
  j_-(x,m) &=& -J_1^0(x,m) = \partial_x \nonumber
 \eea
which represent the Lie algebra sl(2) \cite{tur1}. 
\par The representation of $gl(M+1)$  defined above is irreducible
and has dimension $C^{M+m}_M$, corresponding to the Young
diagram with one line of $M$ elementary boxes.
Series of reducible representations of $gl(M+1)$ can 
trivially be obtained  by considering  the vector space
\be
P(m_1,M)\oplus P(m_2,M)\oplus \cdots \oplus P(m_k,M)
\label{sumgen}
\ee
and building the suitable direct sums of the operators above
\be
   {\rm diag}( J_a^b(\vec x,m_1), J_a^b(\vec x,m_2),\cdots , J_a^b(\vec x,m_k))
+ \delta_a^b {\rm diag}(\gamma_1, \gamma_2, \cdots, \gamma_k)
\label{gl2}
\ee
where we separate explicitely the parameters $\gamma$ entering through
(\ref{j}).
\subsection{Products of representations}
Let  $x_1,\dots,x_M$ and $y_1,\dots, y_N$ denote two sets of 
independent variables, we conveniently define
$P(m,M;n,N)$ as the space of polynomials of overal degree at most
$m$ (resp. $n$) in the variables $x_a$ (resp. $y_b$).
\par Another way to represent $gl(M+1)$ in terms of differential 
operators is to consider the generators corresponding to the 
product of two representations (\ref{j}); 
that is to say,
\be
J^b_a(\vec x,m; \vec y,n) = J^b_a(\vec x,m) + J^b_a(\vec y,n)
\label{sl2d}
\ee
where $\vec x$ and $\vec y$ represent two sets of $M$
independent variables.
Clearly, the operators (\ref{sl2d})
obey the commutation relations (\ref{comj}). They act
on the vector space  $P(m,M;n,M)$. Their action is, however, 
not irreducible.  For shortness, we discuss this statement 
in the case $M=1$, $P(m,1;n,1)$ is abbreviated $P(m;n)$.
\par
One can show easily that the operators (\ref{sl2d})
preserve irreducibly the subspace of $P(m;n)$ defined by
\be
M(m;n) = {\rm Span} \left\lbrace 
({\partial\over{\partial x}} + {\partial\over {\partial y}})^k
x^m y^n \quad , \quad 0 \leq k \leq m+n\right\rbrace
\label{sousespace}
\ee
This vector space is the eigenspace of the Casimir operator
corresponding to the representation of highest spin in the 
decomposition of $P(m;n)$ in subspaces irreducible with 
respect to the action  of (\ref{sl2d}).
The other irreducible representations
result as similar structures with different values of $m$ and $n$;
it is therefore sufficient to deal with (\ref{sousespace}).
\par
Alternatively, $M(m;n)$  can be seen as the kernel  of the operator
\be
K = (x-y) {\partial\over{\partial x}} {\partial \over{\partial y}}
+ n{\partial \over{\partial x}} - m {\partial \over{\partial y}}
\ee
acting on $P(m;n)$.
Remark that $M(m;0)$  (resp. $M(0;n)$) 
is isomorphic to $P(m)$ (resp. $P(n)$);
in that case, the one dimensional operators (\ref{sl2})
are recovered from (\ref{sl2d})  by ignoring
the partial derivative $\partial \over \partial y$ (resp.  
$\partial \over \partial x$).
\subsection{QES operators in two variables}
The QES operators in two real variables are classified in Ref. \cite{gko2};
the authors summarize the seven possible hidden algebras 
in their table 7.
The operators labelled (1.4), (1.10) and (2.3) in this table
are studied independently in Refs. \cite{turshif}, \cite{tur3}. 
The operators labelled (1.1) in the table appear to be new; 
in particular they lead to the only case
for which the invariant module is not 
manifestly a space of polynomials in the two variables.  
In the following we show that the 
formulation (1.1) can be simplified and related
to the algebra (\ref{sl2d}) by means of a suitable change of function.
\par
The operators (1.1) in table 7 of ref.\cite{gko2} read
\begin{eqnarray}
\tilde j_- &=& {\partial\over{\partial x}} + {\partial \over{\partial y}},
\nonumber \\
\tilde j_0 &=& x{\partial\over{\partial x}} + y {\partial\over{\partial y}},
\nonumber \\
\tilde j_+ &=& x^2 {\partial\over{\partial x}} + y^2{\partial\over{\partial y}}
+ {n\over 2} (x-y)
\label{olver11}
\end{eqnarray}
They preserve the space $\tilde M(m;n)$ defined as
\be
\tilde M(m;n) = {\rm span} \left\lbrace (x-y)^{m+{n\over 2}-k} R_k^{m,n}
\left({x+y\over {x-y}}\right)\quad , \quad 0 \leq k \leq 2m+n\right\rbrace
\ee
\be
R_k^{m,n}(t) = {d^k\over{dt^k}} (t-1)^{m+n}(t+1)^m
\ee
Our observation is summarized by the following two formulas
\be
(x-y)^{m+{n\over 2}} \  \tilde j_{\epsilon} \ (x-y)^{-m-{n\over 2}} = 
j_{\epsilon}(x,m) + j_{\epsilon}(y,m+n)
\ee
\be
(x-y)^{m+{n\over 2}} \tilde M(m;n) = M(m;m+n)
\ee
In other words the algebra (\ref{olver11}) is equivalent to the
algebra (\ref{sl2d}) (up to a suitable redefinition of $n$
into $m+n$).
\par The advantage of the new formulation is twofold. 
First the relevant operators form an sl(2) 
diagonal subalgebra of the sl(2)$\times$ sl(2) algebra generated by 
\be
j_{\epsilon}(x,m)\quad , \quad j_{\epsilon}(y,n)\quad , \quad \epsilon =
\pm, 0
\ee
In this respect the form (\ref{sl2d}) of the operators (\ref{olver11})
is clearly related to the fundamental operators (\ref{sl2})
and  is easy to generalize to the case of $M$ variables.  
Second, the vector space of the representation, 
i.e. $M(m;n)$, is a space of polynomials like in
all the other cases of the classification of Ref.\cite{gko2}.
In the next sections, we discuss some possible extensions 
of the algebra (\ref{sl2d}) into graded algebras. The
corresponding operators  are related to  systems of QES  equations. 
\section{Superalgebra spl(M+1,N+1)}
The bosonic part of the superalgebra spl(M+1,N+1) contains the Lie
algebra sl(M+1)$\times$ sl(N+1) \cite{hurni}. 
The fermionic generators
split into two multiplets  ($M+1,\overline{N+1})$ and
($\overline{M+1},N+1)$ under the adjoint action of the
 sl(M+1)$\times$ sl(N+1) subalgebra.
\par
Following the results of \cite{hurni}, 
it appears that the typical representations of spl(M+1,N+1) 
can be constructed by applying suitable combinations of the
fermionic operators to a given irreducible representation of 
sl(M+1)$\times$ sl(N+1), the so called ground floor.
Then the whole representation is generated by applying
the monomials in the fermionic on the ground floor.
The set of vectors attained by means of a monomial of degree
$j$ defines the $j^{th}$ floor.
The anti-commutation relations between the fermionic operators
guarantee that the total number of floors is finite.
\par A natural vector space which can serve
as a basis of representation of  Lie algebra sl(M+1)$\times$ sl(N+1) 
is $P(m,M;n,N)$, defined above.
Then, some irreducible representations of spl(M+1,N+1) can be
constructed in terms of differential operators acting on a direct sum 
\be
\oplus^K_{k=1} P(m_k,M;n_k,N)
\ee
The bosonic generators are the direct sums
\be
    {\rm diag} (J_a^b(\vec x,m_1),\cdots,J_a^b(\vec x,m_K)) \ \ \quad , \ \  
     {\rm diag} (J_a^b(\vec y,n_1),\cdots,J_a^b(\vec y,n_K))
 \ee
 the fermionic ones are built with the tensorial
 operators obtained in \cite{bgk}. They involve in general
  a large  number of 
parameters to be fixed by imposing the commutation relations.
The determination of the possible values for $m_i,m_j,K$ and the explicit
construction of the fermionic generators for generic values of $M,N$ is
a complicated task.
We therefore limit our study to the
particular cases spl(2,1) and spl(2,2).
Our aim is to identify the hidden symmetries 
of quasi exactly solvable systems; 
the simplest of them are related to 
the operators involving a rather low number
of variables and of polynomial components. 
\subsection{The superalgebra spl(2,1)}
The bosonic part of spl(2,1), which is equivalent 
to the (perhaps better known) algebra osp(2,2), 
is the Lie algebra sl(2)$\times$ u(1).  
The structure constant appear rather simple when we
label the bosonic generators $J^b_a$ and the fermionic ones 
$Q_a$ and $\overline Q^a\ (a,b=1,2)$ :
\begin{eqnarray}
&[&J^b_a,J^d_c] = \delta^d_a J^b_c - \delta^b_c J^d_a\\
&[&J^b_a,Q_c] = -\delta^b_c Q_a + \delta^b_a Q_c\\
&[&J^b_a, \overline Q^c] = 
\delta^c_a \overline Q^b - \delta^b_a \overline Q^c\\
&\lbrace& Q_a,\overline Q^b\rbrace = J^b_a\\
&\lbrace& Q_a,Q_b\rbrace = \lbrace \overline Q^a, \overline Q^b\rbrace = 0
\end{eqnarray}
\subsubsection{Representation by differential operators of one variable}
The typical  irreducible representations of spl(2,1) \cite{snr} 
can be expressed  as 
4$\times 4$ matrix differential operators acting 
on the  4-tuple of polynomials of one variable
whose decomposition in floors reads
\begin{eqnarray}
{\rm floor \ 0} &\ & P(m)\\
{\rm floor \ 1} &\ & P(m+1)\oplus P(m-1)\\
{\rm floor \ 2} &\ & P(m)
\end{eqnarray}
\par The operators
\be
\label{qqbar}
q_a(x) = (1,x) \quad ; \quad \overline q_a(x,m) = 
({d\over {dx}}, x{d\over{dx}}-m)
\ee
which naturally connect $P(n)$ with $P(n \pm 1)$
play a crucial role in the construction of the representation.
Using the symplectic metric $\epsilon^{ab}$ 
with $\epsilon^{01}=1$ to raise and lower the
indices, the fermionic generators can be set in the form
\be
Q_a = \left(\begin{array}{cccc}
0&0&0&\\
q_a&0&0&0\\
\overline q_a(n)&0&0&0\\
0&\overline q_a(n+1)&-q_a&0
\end{array}\right)
\ee
\be
\overline Q^a = \left(\begin{array}{cccc}
0&\alpha \overline q^a(n+1)&\beta q^a &0\\
0&0&0&\alpha q^a\\
0&0&0&-\beta \overline q^a(n)\\
0&0&0&0
\end{array}\right)
\ee
with
\be
\alpha = {n-t\over {2(n+1)}}\quad , \beta = {n+2+t\over {2(n+1)}}
\ee
Where $t$ is an arbitrary parameter and we made use of the similarity
transformation to set $Q_a$ in a particularly simple form. The bosonic
operators read then
\be
J^1_0 = -{\rm diag}(j_-(x),j_-(x),j_-(x),j_-(x))
\ee
\be
J^0_1 = {\rm diag} (j_+(x,m),j_+(x,m+1),j_+(x,m-1),j_+(x,m))
\ee
\be
{1\over 2}(J^0_0+J^1_1) = {1\over 2} {\rm diag} (t,t+1,t+2,t+2)
\ee
\be
{1\over 2}(J^0_0 - J^1_1) = 
{\rm diag} (j_0(x,m),j_0(x,m+1),j_0(x,m-1),j_0(x,m))
\ee
The irreducibility of the representation (for generic value of $t$)
and the use of Burnside theorem guarantee that the eight operators
above generate all linear, differential operators preserving
the space $P(m)\oplus P(M+1)\oplus P(M-1) \oplus P(m)$.
\par In the limit $t=-(m+2)$ (i.e. $\beta=0$), 
the upper left $2\times 2$ blocks of the
operators above act invariantly on the space $P(m)\oplus P(m+1)$ 
and lead to a series of  atypical representations.
This formulation of the  Lie superalgebra  osp(2,2)
by differential operators (of one variable) 
preserving $P(m)\oplus P(m+1)$
was first noticed in  \cite{turshif}.
\subsubsection{Representations by differential 
operators of two variables}
By using the formalism of sect.2.1 one can also build
representations of spl(2,1)  in terms
of differential operators of two (or more) variables.              
Let us first discuss the representation by operators
preserving the direct sum
\be
M(m;n) \oplus M(m+1;n) \ \ .
\label{sum1}
\ee
The diagonal generators are of the form (see(\ref{gl2}))
\be
    {\rm diag}(J_a^b(x,m)+J_a^b(y,n) , J_a^b(x,m+1)+J_a^b(y,n))
-  \delta_a^b {\rm diag}(1,0)
\label{jj2}
\ee
(see (\ref{sl2d})), we have to construct the off diagonal ones
which connect the two subspaces of the sum (\ref{sum1}). 
In this case, the relevant tensorial operators  read
\begin{eqnarray}
q_0(x,m;y,n) &=&
  {1\over{m+1}} (m+n+1+(x-y){\partial\over{\partial y}}) \nonumber \\
q_1(x,m;y,n) &=& 
{1\over{m+1}} \left((m+1)  x+ny+y(x-y){\partial
\over{\partial y}}\right)  \nonumber \\
\label{q} \\
\overline q_0(x,m) &=& {\partial\over{\partial x}} \nonumber \\
\overline q_1(x,m) &=& (x {\partial\over{\partial x}} - m )
\ . \label{qqbar2}
\end{eqnarray}
They generalize (\ref{qqbar}) which are recovered 
by setting $n=0$ and dropping
all derivatives ${\partial \over \partial y}$. 
We then obtain the desired representation in terms of
(\ref{jj2}), supplemented by the fermionic generators
\be
Q_{a} = q_{a}(x,m;y,n) \sigma_- \quad  , 
\overline Q^{a} = \overline q^{a}(x,m) \sigma_+ \qquad , \quad a=0,1
\label{op3}
\ee
It has dimension  $2m+2n+3$; again Burnside theorem guarantees
that all differential, linear operators
preserving (\ref{sum1}) are the elements of the envelopping algebra
generated by the eight generators. 
\par
The generic typical representations of spl(2,2) 
can also be formulated in terms 
of the two variables differential operators.
They can be constructed equally well on the vector spaces 
\be
V_I = M_{m,n} \oplus M_{m+1,n}\oplus M_{m-1,n}\oplus M_{m,n}
\ee
\be
V_{II} = M_{m,n} \oplus M_{m+1,n}\oplus M_{m,n+1}\oplus M_{m,n}
\ee
The associated representations are equivalent 
(in agreement with \cite{snr}) but 
the expressions of the generators 
in terms of the partial derivatives is quite different. 
The proof of the equivalence between the representations acting on $V_I$ and 
on $V_{II}$ relies on the fact that the operators 
(\ref{jj2}),(\ref{qqbar2}) are invariant
under the double substitution 
$m \leftrightarrow n$ and $x \leftrightarrow y$.
Therefore the operators $q_a(y,n;x,m)$ 
(resp. $\overline q^a(y,n)$) behave exactly as $q_a(x,m;y,n)$
(resp.  $\overline q^a(x,m)$) under the adjoint action of 
gl(2) algebra (\ref{jj2}).
\subsection{The superalgebra spl(2,2)}
The bosonic part of the superalgebra spl(2,2) is 
sl(2)$\times$ sl(2)$\times$ u(1) 
\cite{hurni}. We note $J^b_a, \tilde J^b_a$ the two sets of
generators of the two sl(2) 
($\sum_a J^a_a = \sum_a \tilde J^a_a = 0$) and
$Y$ the generator related to $u(1)$. Denoting the fermionic generators
by $Q^b_a, \overline Q^b_a$ the commutator relations read as follow
\begin{eqnarray}
&[&J^b_a,J^d_c] = \delta^d_aJ^b_c-\delta^b_e J^d_a\\
&[&\tilde J^b_a,\tilde J^d_c] = \delta^d_a 
\tilde J^b_c-\delta^b_c \tilde J^d_a\\
&[&J^b_a,\tilde J^d_c] = [J^b_a,Y] = [\tilde J^b_a,Y]=0\\
&[&J^b_a,Q^d_c] = -\delta^b_c Q^d_a + {1\over 2} \delta^b_a Q^d_c\\
&[& J^b_a,\overline Q^d_c] = \delta^d_a \overline Q^b_c - {1\over 2} 
\delta^b_a \overline Q^d_c\\
&[&\tilde J_a^b,Q^d_c] = \delta^d_a Q^b_c - {1\over 2} \delta^b_a Q^d_c\\
&[&\tilde J^b_a, \overline Q_c^d] = -\delta^b_c \overline Q^d_a + 
{1\over 2} \delta^b_a \overline Q^d_c\\
&\lbrace& Q^b_a,\overline Q^d_c\rbrace = \delta^b_c J^d_a+\delta^d_a J^d_c
+ {1\over 2} \delta^d_a \delta^b_c Y \\
&\lbrace& Q^b_a, Q^d_c\rbrace = \lbrace \overline Q^b_a,
\overline Q^d_c\rbrace = 0
\end{eqnarray}
\par
The structure of the anticommutators $\lbrace Q,Q \rbrace$ is such
that the generic representation consist of five floors \cite{hurni}.
In  term of polynomial space of two variables $P(m;n)$ 
the representation is organized as follow
\begin{eqnarray}
{\rm floor \ 0} &\ & P(m;n)  \nonumber  \\
{\rm floor \ 1} &\ & P(m+1;n+1)\oplus P(m+1;n-1) \nonumber \\
&\ &\oplus P(m-1;n+1) \oplus P(m-1;n-1) \nonumber \\
{\rm floor \ 2} &\ & P(m+1;n)\oplus P(m;n+2) \oplus P(m;n) \nonumber \\
                &\ &\oplus P(m;n) \oplus P(m;n-2)\oplus P(m-2;n) \nonumber \\
{\rm floor \ 3} &\ & P(m+1;n+1)\oplus P(m+1;n-1) \nonumber \\
&\ &\oplus P(m-1;n+1) \oplus P(m-1;n-1) \nonumber \\
{\rm floor \ 4} &\ & P(m;n)
\end{eqnarray}
The construction of the fermionic generators involves the
determination of 112 parameters. Like for the case spl(2,1) we 
focused our attention on the representations composed out of a 
lower number of levels.
If we restrict to two levels only, we find \cite{hurni}
that the only possibility is the space
$P(1;0)\oplus P(0;1)$, correspondingly $Y=1$.
It was already observed in \cite{bgk} that the 
set of operators preserving the
vector space $P(m;n)\oplus P(k;l)$ do not represent a linear algebra. 
Considering then the representations with three levels, we found a
series of atypical representations with the following dimensions 
\begin{eqnarray}
{\rm floor \ 0} &\ & P(n,n) \nonumber \\
{\rm floor \ 1} &\ & P(n+1,n)\oplus P(n-1,n+1) \nonumber \\
{\rm floor \ 2} &\ & P(n,n)
\end{eqnarray}
The bosonic operators read as follows
\begin{eqnarray}
J^1_0 &=& -{\rm diag} (J_-(x),J_-(x),J_-(x), J_-(x)) \nonumber \\
J^0_0 &=& {\rm diag} (J_0(x,n),J_0(x,n+1),J_0(x,n-1),J_0(x,n) \nonumber \\
J^0_1 &=& {\rm diag} (J_+(x,n), J_+(x,n+1),J_+(x,n-1),J_+(x,n))\nonumber \\
\tilde J^1_0 &=& -{\rm diag} (J_-(y),J_-(y),J_-(y),J_-(y)) \nonumber \\
\tilde J^0_0 &=& {\rm diag} (J_0(y,n),J_0(y,n-1),J_0(y,n+1),J_0(y,n))
\nonumber \\
\tilde J^1_0 &=& {\rm diag} (J_+(y,n),J_+(y,n-1),J_+(y,n+1),J_+(y,n))
\end{eqnarray}
and the value $Y$ is just vanishing for the representation under consideration.
Finally, the fermionic operators are of the form
\be
     Q^b_a = \epsilon^{bc}      V_{ac}\quad, \quad 
\overline Q^b_a = \epsilon^{bc} \overline V_{ca}
\ee
with
\be
V_{ac} = \left(\begin{array}{cccc}
0&0&0&0\\
q_a\overline q_c(n) &0&0&0\\
\overline q_a(n)q_c&0&0&0\\
0&\overline q_a(n+1)q_c&-q_a\overline q_c(n+1)&0
\end{array}\right)
\ee
\be
\overline V_{ac} = {1\over {n+1}} \left(\begin{array}{cccc}
0&-\overline q_a(n+1)q_c &q_a \overline q_c(n) &0\\
0&0&0&q_a\overline q_c(n)\\
0&0&0&\overline q_a(n) q_c\\
0&0&0&0
\end{array}\right)
\ee
where, for shortness, we dropped the variable $x$ (resp. $y$) for the
operator $q_a,\overline q_a$ (resp. $q_c,\overline q_c$). Again, a suitable
similarity transformation was used to set $V_{ac}$ in a form a simple as
possible.
\section{More graded algebras}
The four operators (\ref{qqbar}) play a crucial role
in the construction of the operators preserving the 
spaces $P(m)\oplus P(m+\Delta)$ \cite{bk,bn} and therefore
in the classification of 2$\times$2 QES systems.
It is known that the underlying algebra is non linear 
if $\Delta > 1$.
In this section
we construct the operators preserving the vector space
\be
 M(m;n) \oplus M(m+\Delta;n + \Delta') \quad \ \ . \ \ 
\ee
The relevant diagonal operators are conveniently  chosen according 
to (\ref{gl2})
\bea
&\ &{\rm diag}(J_a^b(x,m)+J_a^b(y,n) ,J_a^b(x,m+\Delta)+J_a^b(y,n+\Delta')) 
\nonumber \\
&- &{1\over 2} \delta_a^b {\rm diag}(1+\Delta+\Delta',1-\Delta-\Delta')
\eea
The construction of the off diagonal ones depends on the relative signs
of $\Delta$ and $\Delta'$.
\subsection{case $\Delta,\Delta' \geq 0$}
The operators connecting the vector space
$M(m;n)$ with $M(m+\Delta;n+\Delta')$ (and vice versa)
can be formulated in terms of products
of  operators (\ref{qqbar2}) where the indices 
$m$ and $n$ in the different factors are appropriately shifted .
The following identities allows one to deal with the 
ambiguities of ordering of
the different factors:
\begin{eqnarray}
&q_b(x,m+1;y,n) \ q_a(x,m;y,n) &=q_a(x,m+1;y,n) q_b(x,m;y,n) \nonumber \\
&q_b(y,n;x,m+1) \ q_a(x,m;y,n) &=q_a(y,n;x,m+1) q_b(x,m;y,n) \nonumber \\
&q_b(y,n;x,m+1) \ q_a(x,m;y,n) &= q_a(x,m;y,n+1) q_b(y,n;x,n) 
\label{iden}
\end{eqnarray}
One can use these identities in order to define the operators 
\be
q(x,[a_k];y,[b_l]) \equiv 
\prod_{l=1}^{\Delta '} 
\prod_{k=1}^{\Delta}
q_{b_l}(y,n+l-1;x,m+ \Delta) 
q_{a_k}(x,m+k-1,n)
\label{qprod}
\ee
symmetrically  in the multi indices 
\bea
&[a_k] \equiv (a_1, \ldots , a_{\Delta '}) \ \ \ , \ \ \ &a_i = 0,1
\nonumber \\
&[b_l] \equiv (b_1, \ldots ,  b_{\Delta }) \ \ \ , \ \ \ &b_i = 0,1
\eea
The operators (\ref{qprod}) connect $M(m;n)$ with $M(m+\Delta;n+\Delta')$
and
\be
Q([a_k],[b_l]) =  q(x,[a_k];y,[b_l]) \sigma_- 
\ee
are the counterparts of the  operators (\ref{q}).
Using the remarks made at the end of the previous
section, one can see that the operators (\ref{qprod})
transform according to the representation of spin $\Delta + \Delta'$
under the adjoint action of the gl(2) represented via (\ref{gl2}).
\par
The operators $\overline Q([a_k],[b_l])$, proportional to
$\sigma_+$,  can be constructed 
in exactly the similar way as for the $Q([a_k],[b_l])$.
Identities like (\ref{iden}) exist among the $ \overline q_a$.
The complete algebra (which is non linear)
 can be obtained following the same lines
as in \cite{bn}.
The anticommutator $\lbrace Q, \overline Q \rbrace$ 
close as polynomial expressions of 
(\ref{gl2}) provided the constants $\gamma$ are choosen such that
\be
J^a_a= {\rm diag} (m+n+\Delta+\Delta'+1 , m+n+1)
\label{top}
\ee
\subsection{The case M(m;n) $\oplus$ M(m+1;n-1)}
If we consider $\Delta$ and $\Delta'$ of opposite signs,
the operators preserving (\ref{sumgen}) do not represent
a Lie super algebra even for 
$\vert \Delta \vert = \vert \Delta' \vert =1$
We studied the operators preserving the vector space 
$M(m;n) \oplus M(m+1;n-1)$
and observed that the underlying algebric stucture is different from those 
obtained in \cite{bn}. Again,
 $J_{a}^{b}(x,m,m+1;y,n,n-1)$ can be used
as a starting point. 
We find it convenient to set $T=0$ in (\ref{gl2})
and add separately the  grading operator 
$T\equiv \sigma_3$.
As far as the off diagonal operators are concerned, we choose
\bea
     &R_{a}^{b} &= 
\overline q_{a}(y,n-1,x,m+1) q^{b}(x,m,y,n) \sigma_-
\nonumber \\
     &\overline R_{a}^{b} &= 
\overline q_{a}(x,m,y,n) q^{b}(y,n-1,x,m+1) \sigma_+
\label{rrbar}
\eea
These tensors are not irreducible with respect to the adjoint action
of the $J$ generators. 
The traces
\be
R_{a}^{a} = 
{m+n+2 \over m+1} ((y-x) {\partial \over \partial y} -n) \ \ \ , \ \ \ 
\overline R_{a}^{a} = 
{m+n+2 \over n+1} ((x-y) {\partial \over \partial x} -m)
\ee
are operators which intertwine the equivalent representations
carried by the spaces $M(m;n)$ and $M(m+1;n-1)$;
they commute with the four operators (\ref{gl2}). 
\par
The  order of the factors $q$ and $\overline q$ entering in
$R$ and $\overline R$ can be reversed by using the identity
\bea
q_{a}(y,n-1,x,m+1) q^{b}(x,m,y,n) &-&
q^{b}(x,m,y,n-1)q_{a}(y,n-1,x,m) \nonumber \\
&=& {1 \over m+n+2} R_{a}^{a}
\eea
The generators $J_{a}^{b},R_{a}^{b}, \overline J_{a}^{b}$
obey   the following commutation 
and anticommutation relations 
\bea
&[ J_{a}^{b} , R_{c}^{d}] &= 
  \delta_{a}^{d} R_{c}^{b} 
- \delta_{c}^{b} R_{a}^{d} \ \ , \nonumber \\
&[ J_{a}^{b} , \overline R_{c}^{d}] &= 
  \delta_{a}^{d} \overline R_{c}^{b} 
- \delta_{c}^{b} \overline R_{a}^{d}
\eea
\bea
\{ R_{a}^{b} ,  \overline R_{c}^{d} \}
&=& {1\over 2} \{ J_{a}^{d} , J_{c}^{b} \}  
+ {T\over 2} (   \delta_{a}^{d}  J_{c}^{b} 
                -  \delta_{c}^{b} J_{a}^{d} ) \nonumber \\ 
&-& {1\over 2} (   \delta_{a}^{b}  J_{c}^{d} 
                +  \delta_{c}^{d} J_{a}^{b} ) 
- {1\over 2}  \delta_{a}^{b}  \delta_{c}^{d}   
\label{corrbar}
\eea
These  relations do not define an abstract algebra.
In order to fulfil the associativity conditions
(i.e. the Jacobi identities), the  anticommutators between 
two $R$ (two $\overline R$), which vanish for (\ref{rrbar}) 
have to be implemented in a non trivial way \cite{bk},\cite{bn};
equation (\ref{corrbar}) indicates that the underlying algebra 
is non linear.
\section{Deformation of spl(2,1)}
If we want to describe in abstract terms the
algebraic structure  underlying the QES finite difference equations,
some deformations of the algebras discussed above seem to emerge
in a natural way.
For scalar equations the relevant deformation was pointed out
some time ago \cite {tur2}; it is related to the 
``Witten type II'' deformation  of sl(2) \cite{witt,zachos}
Since the most relevant examples of QES systems are 
related to spl(2,1), it is natural to try to construct deformations
of spl(2,1) whose representations can be formulated in terms of 
finite difference operators.
\par
One deformation of spl(2,1) was obtained in \cite{deguchi} 
(for more general graded algebra see \cite{vinet});
This deformation is such that the commutators of some generators 
close into transcendental functions of other generators.
In a previous paper \cite{bgk} some representations
of this deformation of spl(2,1) 
were formulated in terms of  finite difference operators.
\par
More recently \cite{bh} a family of quommutator 
deformations of the superalgebra
spl(p,1) was constructed.
The quommutators and anti-quommutators are defined respectively
as
\be
   [A , B]_q = AB - q BA \ \ \ \ , \ \ \ \ \{A , B \}_q = AB + qBA
\ee
where $q$ is the deformation parameter.
The quommutator deformations are closer, in the spirit,
to the Witten type II deformation of sl(2). 
They lead in particular to natural ordering rules on the generators. 
In contrast to the deformation \cite{deguchi}, 
exchanging the order of quommutating generators can be performed
while keeping the polynomial character of the expression.
\par
Here we present with the tensor notations the deformation of spl(2,1) 
which appears to be the most relevant for QES difference equations. 
It is expressed as follows, ($\mu,\nu,\alpha \dots = 0,1$)
\be
  \{ Q_{\mu} , Q_{\nu} \}_{q^{\nu-\mu}} = 0 \ \ \  , \ \ \  
 \{ \overline Q^{\mu} , \overline Q^{\nu} \}_{q^{\nu-\mu}} = 0 
\ee
\be
 \{ Q_{\mu} , \overline Q^{\nu} \} = J_{\mu}^{\nu}  
\label{defj}
\ee
\be
[ J_{\mu}^{\nu} , Q_{\alpha} ]_{q^{\alpha-\nu}} =
q^{{\alpha-\nu-1}\over 2}  
 (\delta_{\mu}^{\nu} Q_{\alpha} - \delta_{\alpha}^{\nu}  Q_{\mu}) 
\ee
\be
[ J_{\mu}^{\nu} , \overline Q^{\alpha} ]_{q^{\mu -\alpha}} =  
q^{{\mu-\alpha-1}\over 2}  
 (\delta^{\alpha}_{\mu}  \overline Q^{\nu} -
    \delta_{\mu}^{\nu} \overline Q^{\alpha}  ) 
\ee
\bea
[ J_{\mu}^{\nu} , J_{\alpha}^{\beta}]_{q^s} &=& 
 q^{{s-r-1}\over 2}(\delta_{\mu}^{\beta} J_{\alpha}^{\nu} 
- q^{r} \delta_{\alpha}^{\nu} J_{\mu}^{\beta}) \nonumber \\
&+& {q-1 \over q^2} (Q_0 \overline Q^0 + Q_1 \overline Q^1)
\delta_{\alpha}^{\nu} \delta_{\mu}^{\beta}
(1-\delta_{\mu}^{\nu} \delta_{\alpha}^{\beta})
\label{gl2def}
\eea
with $s\equiv \nu + \alpha - \mu - \beta$ and
$r \equiv (\nu - \beta)(\mu - \alpha)$.
\par
All Jacobi identities are obeyed by these relations.
The last commutator indicates that 
(like for  the deformation \cite{deguchi})
the bosonic generators do not 
close into a gl(2) subalgebra .
Neglecting all fermionic operators in the above formulas
leads to one deformation of gl(2) found in \cite{fairlie}. 
\par
It is possible to construct two independent expressions,
quadratic in the generators, which q-commute with the generators.
These ``q-Casimir'' operators read 
\bea
   &C_1 = &Q_0 \overline Q^0 + Q_1 \overline Q^1
         + q J_1^0 J_0^1 - J_0^0 J_1^1 - J_0^0 \nonumber \\
   &C_2 = &(q-1)^2 (Q_0 \overline Q^0 + Q_1 \overline Q^1) \nonumber \\
        & &+ q(q-1)^2 J_1^0 J_0^1 + (q-1) J_1^1 -q(q-1)J_0^0 -1 
\eea
and obey the following commutation properties (for i=1,2)
\bea
   &[C_i, J_{\mu}^{\nu} ]_{q^{2(\mu - \nu)}} &= 0 \ \ \ ,  \nonumber \\
   &[C_i, Q_{\mu}]_{q^{2\mu - 1}} &= 0 \ \ \ ,  \nonumber \\
   &[C_i, Q^{\mu}]_{q^{1-2\mu }} &= 0 \ \ \ 
\eea
It follows that any function of the ratio $C_1/C_2$ commute 
with the generators and constitute a conventional Casimir.
\par       
We further constructed the representations of the algebra above
which are relevant for systems of finite difference QES equations.
To describe them we define a finite 
difference operator  $D_q$.  
\be
       D_q f(x) =  { f(x) - f(qx) \over (1-q)x} \ \ \ , \ \ \ 
       D_q x^n = [n]_q x^{n-1} \ \ \ , \ \ \   
       [n]_q \equiv {1 - q^n \over 1-q}  
\ee 
\par
The simplest of these realizations are characterized by a positive
integer $n$ and  act on the vector space
\be
P(n-1) \oplus P(n) \ \ .
\label{pp}
\ee
 Adopting $x$ as the variable,  
the fermionic generators are represented by
\bea
   &Q_0 = q^{-{n\over 2}} \sigma_- \ \ \ , \ \ \  
   &Q_1 = -x \sigma_- \nonumber \\ 
   &\overline Q^0 = q^{-{n\over 2}} (x D_q - [n]_q) \sigma_+ \ \ \ , \ \ \  
   &\overline Q^1 = D_q \sigma_+ 
\label{repp}
\eea
The bosonic operators can be constructed easily from (\ref{defj}) but
we write them for completeness
\be
J_0^0 = q^{-n} ((x D_q - [n]_q) \unop_2 \ \ \ \ , \ \ \ 
J_1^1 = (-1) {\rm diag } ( q x D_q + 1 , x D_q)
\ee
\be
J_0^1 = q^{-{n\over 2}}  D_q \unop_2 \ \ \ \ , \ \ \ 
J_1^0 = (-1)q^{-{n\over 2}}  {\rm diag } 
( qx(x D_q - [n-1]_q) , x(x D_q - [n]_q) )
\ee
The enveloping algebra constructed with the eight generators above
contains all finite difference operators preserving $P(n-1) \oplus P(n)$.
\par
As in the underformed case, there also exist representations
which act on the vector space 
\be
     P(n) \oplus P(n+1) \oplus P(n-1) \oplus P(n)
\label{pppp}
\ee
The fermionic generators are realized as follows
\be
Q_0 =
\left(\begin{array}{cccc}
0     &0     &0   &0\\
1     &0     &0   &0\\
D_q   &0     &0   &0\\
0     &-D_q  &1   &0
\end{array}\right)
\nonumber
\ee
\be
Q_1 =
\left(\begin{array}{cccc}
0     &0     &0   &0\\
x     &0     &0   &0\\
q^{-n}\delta(n)   &0     &0   &0\\
0     &-q^{-n-1}\delta(n+1)  &x   &0
\end{array}\right)
\nonumber
\ee
\be
\overline Q^0 =
\left(\begin{array}{cccc}
0  &\lambda \delta(n+1) &(1-\lambda q^{n+1})x    &0    \\
0  &0                   &0   &(\lambda q^{n+1}-1)x       \\
0  &0                   &0   &q\lambda \delta(n)         \\
0  &0                   &0   &0   
\end{array}\right)
\nonumber
\ee
\be
\overline Q^1 =
\left(\begin{array}{cccc}
0  &\lambda q^{n+1}D_q  &(1-\lambda q^{n+1})      &0   \\
0  &0                   &0   &\lambda q^{n+1} -1       \\
0  &0                   &0   &\lambda q^{n+1} D_q      \\
0  &0                   &0   &0   
\end{array}\right)
\label{repppp}
\ee
where $\lambda$ is an arbitrary complex parameter and, 
for shortness, we used $\delta(n) \equiv xD_q - [n]_q$. 
Once more, the invariance  under the similarity
transformations has been exploited to set $Q_0$ in a form
as simple as possible.
In the limit $\lambda \rightarrow q^{-n-1}$ the representation
(\ref{repppp}) becomes reducible and decomposes 
into atypical ones of the form (\ref{repp}).
\par
The deformation  of spl(2,1) presented above leads to  
simple  normal ordering rules for the generators.
In this respect, it is very appropriate 
for the classification of the finite 
difference operators preserving the space (\ref{pp}) or (\ref{pppp})
and of the corresponding QES systems.
The normal ordering rules 
associated with the deformation used in ref.\cite{bgk} 
are not as transparent since they depend non polynomially
of some operators.
\section{Concluding remarks}
The most interesting examples of QES systems are related to
the algebra spl(2,1), e.g. the relativistic Coulomb problem
and the stability of the sphaleron in the abelian Higgs model
\cite{bk}.
Very recently, the hidden algebras of the supersymmetric
versions of the Calogero and Sutherland models were shown
to be the superalgebras spl(N+1,N) \cite{brink}. 
\par
 Further realizations of these algebras, and more generally
of the superalgebras spl(p,q),
in terms of differential operators therefore desserve some attention.
Here, we present realizations of spl(2,1) formulated in terms of 
differential operators in two variables; extensions to the
case of an arbitrary number of variables can be achieved
in a straightforward way.
We also constructed a series of atypical representations
of spl(2,2) in terms of operators acting on a set of 4-tuple 
of polynomials in two variables.
The labelling used for the generators clearly exhibits their
tensorial structure under their bosonic subalgebra and 
provides very naturally the building blocks for the construction
of series of (non linear) graded algebras preserving more
general vector spaces of polynomials.  
\par
Witten type deformations attracted recently some attention
(see e.g. \cite{chung} for  osp(1,2)).
The deformation of spl(2,1) presented here is of this type
and it admits representations which are directly relevant for
the study of QES finite difference systems.
The correspondance between discrete operators and finite
difference one is discussed in \cite{wieg} where also some
some applications in discrete quantum mechanics are pointed out.
\par
The existence of a coproduct for this type of deformation
would allow to adapt the construction of section 2 to
finite difference equations. This is nice an application
of the coproduct that we plan to address later.

\end{document}